\documentclass[draft]{amsart}
\usepackage{amsmath}
\usepackage{amssymb}
\usepackage{amscd,amssymb,amsmath,latexsym,enumerate}
\usepackage[mathscr]{euscript}
\usepackage{epsfig}
\usepackage{fancybox}
\usepackage{yhmath}
\usepackage{verbatim}
\usepackage{graphicx,bm,units,yfonts}
\usepackage{phonetic}
\usepackage{mathtools}
\theoremstyle{proclaim}
\newtheorem{theorem}{Theorem}[section]

\newtheorem{proposition}[theorem]{Proposition}
\theoremstyle{fancyproclaim}

\theoremstyle{statement}
\newtheorem{remark}[theorem]{Remark}
\newtheorem{definition}[theorem]{Definition}

\theoremstyle{fancystatement}

\numberwithin{equation}{section}

\providecommand{\AMS}{$\mathcal{A}$\kern-.1667em%
\lower.25em\hbox{$\mathcal{M}$}\kern-.125em$\mathcal{S}$}

\begin{document}
%\issueinfo{vv}{n}{yyyy} 

%\commby{Editor}
%\pagespan{1}{23}

%\date{Month dd, yyyy}

%\revision{Month dd, yyyy}

\title[Intrinsic Chern-Connes Characters]{Intrinsic Chern-Connes Characters for Crossed Products by $\mathbb Z^d$} 

\author{Emil Prodan}

\address{Emil Prodan, Department of Physics, Yeshiva University, New York, NY 10016, USA}

\email{prodan@yu.edu}

\begin{abstract} By imbedding $\mathcal A \rtimes_\xi \mathbb Z^d$ into the $C^\ast$-algebra of adjointable operators over the standard Hilbert $\mathcal A$-module $\mathcal H_{\mathcal A}$, we define Fredholm modules and Chern-Connes characters that are intrinsic to the $C^\ast$-dynamical system $(\mathcal A,\xi,\mathbb Z^d)$.  We introduce a $\widehat{\mathcal T}$-index for the generalized Fredholm operators over $\mathcal H_{\mathcal A}$ and develop a generalized Fedosov principle and formula. We prove an index formula for the pairing of the characters with $\bm K_0(\mathcal A \rtimes_\xi \mathbb Z^d)$ and conclude that the pairing is in the image of $\bm K_0(\mathcal A)$ under the trace $\widehat{\mathcal T}$. A local index formula enables new applications in condensed matter physics.
\end{abstract}

\thanks{This work was supported by the U.S. NSF grant DMR-1056168.}

\maketitle

\section*{INTRODUCTION}

Consider a classical dynamical system $(\Omega, \xi)$, where $\xi=(\xi_1,\ldots,\xi_d)$ is a system of $d$-commuting homeomorphisms ($d =$ even). Let $\big ( C(\Omega), \xi,\mathbb Z^d\big)$ be its dual $C^\ast$-dynamical system and $C(\Omega) \rtimes_\xi \mathbb Z^d$ the canonically associated crossed-product. Let $\pi_\omega$ be the standard representation of the crossed product on $\ell^2(\mathbb Z^d)$:
\begin{equation}\label{StandardRep}
\begin{array}{l}
\Big(\pi_\omega \Big (\sum_{\bm q} \phi_{\bm q} \cdot \bm q \Big) \psi \Big)_{\bm x} = \sum_{\bm q} \phi_{\bm q}(\xi_{\bm x} \omega) \psi_{\bm x - \bm q},
\end{array}
\end{equation}
and $D_{\bm x_0}={\bm \gamma}\otimes (\bm X+\bm x_0)$ be the (shifted-) Dirac operator on $\mathbb C^{2^{\frac{d}{2}}} \otimes \ell^2(\mathbb Z^d)$. Then $\big(\mathbb C^{2^\frac{d}{2}} \otimes \ell^2(\mathbb Z^d), \pi_\omega, F_{\bm x_0},\gamma\big)$ with $F_{\bm x_0} = \mathrm{sign}(D_{\bm x_0})$ is a natural even Fredholm module over $C(\Omega) \rtimes_\xi \mathbb Z^d$ and its Chern-Connes character, defined as the cohomology class of the cyclic cocycle:
\begin{equation}
 \tau_d(\bm a_0, \ldots \bm a_d) = \nicefrac{1}{2} \ \mathrm{Tr}\big \{\gamma F_{\bm x_0} [F_{\bm x_0}, \pi_\omega(\bm a_0)] \ldots [F_{\bm x_0}, \pi_\omega(\bm a_d)] \big \},
 \end{equation}
 pairs well and integrally with the $\bm K_0$-group \cite{CONNES:1985cc}:
\begin{equation}\label{Pairing}
 \bm K_0(C(\Omega) \rtimes_\xi \mathbb Z^d) \ni [\bm p]_0\rightarrow \tau_d(\bm p, \ldots \bm p)=\mathrm{Index} \big\{\pi_\omega^+(\bm p) F_{\bm x_0} \pi_\omega^-(\bm p)\big \} \in \mathbb Z.
\end{equation}
In \cite{BELLISSARD:1994xj,ProdanJPA2013hg}, under certain \emph{optimal conditions}, the local formula:
\begin{equation}\label{L}
\tau_d(\bm a_0, \ldots \bm a_d) = \Lambda_d\sum_{\rho \in S_d} (-1)^{\rho} \mathcal T\{\bm a_0 \partial_{\rho_1} \bm a_1 \ldots \partial_{\rho_d} \bm a_d\}, \ \Big(\Lambda_d=\frac{(2 \pi \sqrt{-1})^\frac{d}{2}}{(d/2) !}\Big)
\end{equation}
was proven by \emph{elementary means}. Above, $S_d$ is the group of permutations and $(\bm \partial, \mathcal T)$ is the noncommutative differential calculus over $C(\Omega) \rtimes_\xi \mathbb Z^d$.

The identity \ref{CentralId} played a key role in this computation. In $d=2$, this identity is due to Connes \cite{CONNES:1985cc}, who used it to compute the 2-dimensional Chern characters of the convolution algebras $C_c^\infty(\mathbb R^2)$ and $C_c^\infty(SL(2,\mathbb R))$. The local index formula of \cite{BELLISSARD:1994xj,ProdanJPA2013hg} can be seen as a particularization of the generic Connes-Moscovici formula \cite{ConnesGFA1995re} and its later extensions \cite{HigsonICTP2003gf,CareyAdvMath2006bb,CareyAdvMath2006aa,CareyNCGPhysics,CareyAMS2014jf}. A recent related work is \cite{AnderssonARXIV2014}, where local index formulas for Rieffel deformed crossed products are derived. Among these works, only \cite{CareyAMS2014jf} seems to cover the general settings of \cite{BELLISSARD:1994xj,ProdanJPA2013hg} where the index formulas were shown to hold over certain noncommutative Sobolev spaces. Notice also that \ref{L} gives a local formula for the character itself and not just of the pairing. 

There is some interest from the condensed matter physics community in the results summarized by \ref{StandardRep}-\ref{L} because all non-interacting quantum lattice models of homogeneous materials can be generated as representations of crossed products by $\mathbb Z^d$. The righthand side of \ref{L} relates to the transport coefficients of real materials, in particular, it has certain relevance for topological insulators \cite{ProdanTQM2014bv}. While \cite{BELLISSARD:1994xj,ProdanJPA2013hg} explained some outstanding properties of these materials in the presence of strong disorder and magnetic fields, the approach is limited to the single particle theory of solids, where the electron-electron interaction is treated as a mean-field correction. In fact, even before considering the electron-electron interaction, one needs to address the following shortcomings:
\begin{enumerate}[(a)]
 \item The formalism works only for crossed products of commutative algebras.
\item The $\bm K$-groups of many crossed-products used in condensed matter can be fully resolved by the top and the lower Chern numbers. \cite{ProdanJPA2013hg} produced a non-commutative theory only for the top Chern number.
\item The pairing of the characters with the $\bm K$-groups is always integral, hence not relevant for the sequences of fractional topological phases, such as the fractional Chern insulators \cite{ParameswaranCRP2013gg}.
\end{enumerate}
Removing these deficiencies while maintaining the elementary character of the calculation was the main motivation for the present work.\medskip

The root of the problem is the representation of the Fredholm modules on Hilbert spaces and a natural  cure is provided by $\bm K\bm K$-theory. While exploring this path, we discovered that the crossed-product $\mathcal A \rtimes_\xi \mathbb Z^d$ can be canonically imbedded in the $C^\ast$-algebra of adjointable operators over the standard Hilbert $\mathcal A$-module $\mathcal H_{\mathcal A}$. A Dirac operator over $\mathcal H_{\mathcal A}$ can be naturally defined from the action of $\mathbb Z^d$ after tensoring with an appropriate Clifford algebra. Furthermore, if $\mathcal A$ posses a continuous trace, then this trace can be naturally promoted to a lower semicountinuous trace $\widehat{\mathcal T}$ over the imbedding algebra. The definitions of the intrinsic Fredholm modules, of the notion of summability and of the intrinsic Chern-Connes characters are then fairly straightforward (see Defs.~\ref{GFredholmModule}, \ref{SummaDef} and \ref{ChernConnesCh}, respectively).

One outcome of the approach (see Th.~\ref{KKMap}) is that the operator replacing the one appearing inside the Index in \ref{Pairing} is now a generalized Fredholm operator over $\mathcal H_{\mathcal A}$. As such, Mingo's index \cite{MingoTAMS1987fg} for $C^\ast$-modules provides a $KK$-map from $\bm K_0(\mathcal A \rtimes_\xi \mathbb Z^d)$ to $\bm K_0(\mathcal A)$. One could imagine the possibility of the crossed product $\mathcal A \rtimes_\xi \mathbb Z^d$ being imbedded in the algebra of adjointable operators over $\mathcal H_{\mathcal B}$ of another $C^\ast$-algebra, in which case one will perhaps obtain a $\bm K\bm K$-map into $\bm K_0(\mathcal B)$. This could be a useful tool for the computation of the $\bm K_0$-groups. However, what we find interesting about our construction is that it is intrinsic, in the sense that the entire construction is natural and relying entirely on data from the $C^\ast$-dynamical system $(\mathcal A,\xi,\mathbb Z^d)$. 

In analogy with the Breuer-Fredholm index for von Neumann algebras \cite{BreuerMA1968aa,BreuerMA1968bb} and the later generalizations \cite{OlsenMAMS1984bc}, and especially \cite{PhillipsFIC1997fj}, we define a numerical $\widehat{\mathcal T}$-index for the multiplier algebra, by applying the trace $\widehat{\mathcal T}$ on the Mindex of generalized Fredholm operators over $\mathcal H_{\mathcal A}$. This numerical index takes values in the image of  $\bm K_0(\mathcal A)$ under the trace $\widehat{\mathcal T}$ and the structure of this sub-group of $\mathbb R$ can be far more complex than $\mathbb Z$. For example, if $\mathcal A$ is the irrational rotational algebra $\mathcal A_\theta$, then this image is at least $(\mathbb Z + \theta \mathbb Z)\cap [0,1]$ \cite{RieffelPJM1981bf}. This addresses point (c) above. In order to compute the $\widehat{\mathcal T}$-index, we develop a generalized Fedosov principle and a Fedosov formula (see Th.~\ref{FedosovTh}), which enables us to prove an index formula for the pairing of the intrinsic Chern-Connes characters with $\bm K_0(\mathcal A \rtimes_\xi \mathbb Z^d)$. Hence, this pairing is in the image of $\bm K_0(\mathcal A)$ under the trace $\widehat{\mathcal T}$.

Lastly, using the same elementary methods as in \cite{ProdanJPA2013hg}, we derive a local formula for the Chern-Connes cocycle, similar to \ref{L} (see Th.~\ref{LocalFormulaTh}). The theory now covers the lower Chern numbers, though only in the regime of weak disorder. An application for disordered topological insulators in 3 space-dimensions is provided in the last Chapter, where an index formula for the so called weak topological invariants is derived and predictions about their possible values are made.

\section{PRELIMINARIES}

Let $(\mathcal A, \xi, \mathbb Z^d)$ be a $C^\ast$-dynamical and $\mathcal A \rtimes_\xi \mathbb Z^d$ its canonical crossed-product. The $C^\ast$-algebra $\mathcal A$ is assumed separable and to posses a unit and a faithful, continuous and $\xi$-invariant trace $\mathcal T_{\mathcal A}$. This is the only input we need for the definition and characterization of the intrinsic Chern-Connes characters. 

Throughout our presentation, we follow closely the notation from Davidson's monograph \cite{DavidsonBook1996bv}. In particular, the core algebra $\mathcal A \mathbb Z^d$ will be represented by formal finite sequences:
\begin{equation}
\bm a=\sum_{\bm q} a_{\bm q} \cdot \bm q, \ \ a_{\bm q} \in \mathcal A, \ \bm q \in \mathbb Z^d,
\end{equation}
together with the standard algebraic operations. The algebra $\mathcal A$ is imbedded in the crossed-product as:
\begin{equation}
\mathcal A \ni a \rightarrow \bm a = a \cdot \bm 0 \in \mathcal A \rtimes_\xi \mathbb Z^d,
\end{equation}
and the additive group as:
\begin{equation}
\mathbb Z^d \ni \bm q \rightarrow \bm u_{\bm q} = 1 \cdot \bm q \in \mathcal A \rtimes_\xi \mathbb Z^d.
\end{equation}
The first imbedding is isometric and the second one is in the group of unitaries of $\mathcal A \rtimes_\xi \mathbb Z^d$.

The Fourier calculus is defined by the group of automorphisms $\big\{\rho_{\bm \lambda}\big \}_{\bm \lambda \in \mathbb T^d}$ on $\mathcal A \rtimes_\xi \mathbb Z^d$ \cite{DavidsonBook1996bv}, which act on the core algebra as:
\begin{equation}\label{Rho}
\mathbb T^d \ni \bm \lambda \rightarrow \rho_{\bm \lambda}(\bm a) = \sum_{\bm q} \bm \lambda^{\bm q} a_{\bm q} \cdot \bm q, \ \ \bm \lambda^{\bm q}=\lambda_1^{q_1} \ldots \lambda_d^{q_d}.
\end{equation}
The Fourier coefficients of $\bm a \in \mathcal A \rtimes_\xi \mathbb Z^d$ are defined by the Riemann integral:
\begin{equation}
\Phi_{\bm q}(\bm a) = \int_{\mathbb T^d} d \mu(\bm \lambda) \ \rho_{\bm \lambda} (\bm a \bm u_{\bm q}^{-1}) =\int_{\mathbb T^d} d \mu(\bm \lambda) \ \lambda^{-\bm q}\rho_{\bm \lambda} (\bm a) \bm u_{\bm q}^{-1},
\end{equation}
where $\mu$ is the Haar measure on $\mathbb T^d$. The Fourier coefficients will be seen as elements of $\mathcal A$.

\begin{proposition}[\cite{DavidsonBook1996bv}, pg.~223]\label{Feher} The Ces\`{a}ro sums:
\begin{equation}
\bm a_N=\sum_{q_1=-N}^N \ldots \sum_{q_d=-N}^N \prod_{j=1}^d \left (1-\frac{|q_j|}{N+1} \right ) \Phi_{\bm q}(\bm a) \cdot \bm q
\end{equation}
converge in norm to $\bm a \in \mathcal A \rtimes_\xi \mathbb Z^d$ as $N \rightarrow \infty$.
\end{proposition}
\noindent Hence, generic elements from $\mathcal A \rtimes_\xi \mathbb Z^d$ can be representated as a Fourier series:
\begin{equation}
\bm a=\sum_{\bm q \in \mathbb Z^d} \Phi_{\bm q}(\bm a) \cdot \bm q,
\end{equation}
where the infinite sum must be interpreted via \ref{Feher}.

The Fourier calculus generates a canonical faithful and continuous trace on $\mathcal A \rtimes_\xi \mathbb Z^d$ \cite{DavidsonBook1996bv}:
\begin{equation}\label{Trace0}
\mathcal T \{ \bm a \} = \mathcal T_{\mathcal A}\{\Phi_{\bm 0}(\bm a)\}.
\end{equation}
It also generates a set of un-bounded derivations, $\bm \partial=(\partial_1, \ldots, \partial_d)$, through the generators of the $d$-parameter group of automorphisms $\{\rho_{\bm \lambda}\}_{\bm \lambda \in \mathbb T^d}$. The derivations act as:
\begin{equation}
\bm \partial \bm a =\imath \sum_{\bm q \in \mathbb Z^d} \bm q \Phi_{\bm q}(\bm a) \cdot \bm q, \ (\imath=\sqrt{-1}).
\end{equation}
 Together, $(\bm \partial, \mathcal T)$ define the non-commutative calculus over $\mathcal A \rtimes_\xi \mathbb Z^d$. 

Generically, the cocycles can be defined only on a pre $C^\ast$-sub-algebra of $\mathcal A \rtimes \mathbb Z^d$ \cite{Connes:1994wk}, which can be generated by various means \cite{SchweitzerIJM1993re,RennieKTh2003vj}. Below we describe one such sub-algebra which we find particularly convenient for the calculations to follow.

\begin{proposition} Consider the set $(\mathcal A \rtimes_\xi \mathbb Z^d)_{\mathrm{loc}}$:
\begin{equation}\label{ALoc1}
 \Big \{\bm a \in \mathcal A \rtimes_\xi \mathbb Z^d \ | \ \exists \ \alpha<1, \ s.t.  \ \sup_{\bm q} \big (\alpha^{-|\bm q|} \| \Phi_{\bm q}(\bm a)\| \big ) <  \infty \Big \}. 
\end{equation}
Then:
\begin{enumerate}
\item The set $(\mathcal A \rtimes_\xi \mathbb Z^d)_{\mathrm{loc}}$ can be equivalently characterized as:
\begin{equation}\label{ALoc2}
\Big \{\bm a \in \mathcal A \rtimes_\xi \mathbb Z^d \ | \ \rho_{\bm \lambda}(\bm a) \ \mbox{analytic of $\bm \lambda$ in a strip around} \ \mathbb T^d \Big \}. 
\end{equation}
\item  When equipped with the algebraic operations, $(\mathcal A \rtimes_\xi \mathbb Z^d)_{\mathrm{loc}}$ becomes a dense sub-algebra of $\mathcal A \rtimes_\xi \mathbb Z^d$. \medskip
\item This sub-algebra is stable under the holomorphic functional calculus.
\end{enumerate}
\end{proposition}

\proof ({\it i}) Note that $\rho_{\bm \lambda}$ is entire on $\mathcal A \mathbb Z^d$. If $\bm b \in \mathcal A \mathbb Z^d$ and $\bm a \in \mathcal A \rtimes_\xi \mathbb Z^d$ with $\rho_{\bm \alpha}(\bm a)$ analytic of $\bm \alpha$ in a finite strip around $\mathbb T^d$, then $\rho_{\bm \alpha}(\bm a \bm b)$ can be analytically continued in the same strip via $\rho_{\bm \alpha}(\bm a \bm b) = \rho_{\bm \alpha}(\bm a)\rho_{\bm \alpha}(\bm b)$. Likewise, if $\bm \lambda \in \mathbb T^d$, then $\rho_{\bm \alpha \bm \lambda}(\bm a)$ can be analytically continued in $\bm \alpha$ via: 
\begin{equation}
\rho_{\bm \alpha \bm \lambda}(\bm a)=\rho_{\bm \lambda} \circ \rho_{\bm \alpha}(\bm a)=\rho_{\bm \alpha} \circ \rho_{\bm \lambda}(\bm a).
\end{equation}
The last equality holds because the analytic continuations are unique \cite{BlumTAMS1955bv}. Then:
\begin{align}
\Phi_{\bm q}\big (\rho_{\bm \alpha}(\bm a)\big ) & = \int_{\mathbb T^d} d\mu(\bm \lambda) \ \rho_{\bm \lambda} \big (\rho_{\bm \alpha}(\bm a)  \bm u_{-\bm q}\big )  =\alpha^{\bm q} \int_{\mathbb T^d} d\mu(\bm \lambda) \ \rho_{\bm \alpha} \big (\rho_{\bm \lambda}(\bm a  \bm u_{-\bm q})\big ).
\end{align}
We can exchange $\rho_{\bm \alpha}$ and the integral by using the dominated convergence theorem, to conclude:
\begin{equation}
\Phi_{\bm q}\big (\rho_{\bm \alpha}(\bm a)\big ) =\bm \alpha^{\bm q} \rho_{\bm \alpha}\big (\Phi_{\bm q}(\bm a) \big ) = \bm \alpha^{\bm q} \Phi_{\bm q}(\bm a).
\end{equation}
This identity is valid for any $\bm \alpha$ in a strip around $\mathbb T^d$, and this strip is independent of $\bm q$. Taking $\alpha_i = \alpha^{-q_i/|\bm q|}$ with $|\alpha|$ close-enough to 1, we obtain:
\begin{equation}
\|\Phi_{\bm q}(\bm a)\| \leq \alpha^{|\bm q|} \|\rho_{\bm \alpha}(\bm a)\|, \ \forall \bm q \in \mathbb Z^d.
\end{equation}

Now assume $\bm a\in (\mathcal A \rtimes_\xi \mathbb Z^d)_{\mathrm{loc}}$. We can apply $\rho_{\bm \alpha}$ on the Ces\`{a}ro sums:
\begin{equation}
\rho_{\bm \alpha}(\bm a_{N})=\sum_{q_1=-N}^N \ldots \sum_{q_d=-N}^N \prod_{j=1}^d \left (1-\frac{|q_j|}{N+1} \right ) \bm \alpha^{\bm q} \Phi_{\bm q}(\bm a) \cdot \bm q,
\end{equation}
and the righthand-side and its derivatives with respect to $\bm \alpha$ can be seen to be absolutely norm-convergent as $N \rightarrow \infty$, for $\bm \alpha$ in a thin-enough strip around $\mathbb T^d$. The statement then follows.\medskip

({\it ii}) $(\mathcal A \rtimes_\xi \mathbb Z^d)_{\mathrm{loc}}$ is closed under the algebraic operations, which can be seen directly from representation \ref{ALoc2}. Indeed, if $\rho_{\bm \alpha}(\bm a)$ and $\rho_{\bm \alpha}(\bm b)$ are analytic in a strip around $\mathbb T^d$, then $\rho_{\bm \alpha}(\bm a \bm b)$ can be analytically continued over the same strip via $\rho_{\bm \alpha}(\bm a \bm b)=\rho_{\bm \alpha}(\bm a) \rho_{\bm \alpha}(\bm b)$. $(\mathcal A \rtimes_\xi \mathbb Z^d)_{\mathrm{loc}}$ is dense in $\mathcal A \rtimes_\xi \mathbb Z^d$ because it contains $\mathcal A \mathbb Z^d$.\medskip

({\it iii}) Consider $ \bm a \in (\mathcal A \rtimes_\xi \mathbb Z^d)_{\mathrm{loc}}$ which is invertible in $\mathcal A \rtimes_\xi \mathbb Z^d$. We need to show that its inverse belongs to  $(\mathcal A \rtimes_\xi \mathbb Z^d)_{\mathrm{loc}}$. Since the latter is dense in $\mathcal A \rtimes_\xi \mathbb Z^d$, there is $\bm b$ from $(\mathcal A \rtimes_\xi \mathbb Z^d)_{\mathrm{loc}}$ such that $\|\bm a \bm b - 1\|<1$. Then $\bm a \bm b=1-(1-\bm a \bm b)$ is invertible, hence $\bm b$ is invertible and $\bm a^{-1}=\bm b(\bm a \bm b)^{-1}$. Taking $\bm r =1-\bm a \bm b$, the problem is reduced to showing that $(1-\bm r)^{-1}$ belongs to $(\mathcal A \rtimes_\xi \mathbb Z^d)_{\mathrm{loc}}$.  Recall that $\|\bm r\|<1$ and $\bm r\in (\mathcal A \rtimes_\xi \mathbb Z^d)_{\mathrm{loc}}$. Now, $\rho_{\bm \alpha}(\bm r)$ is analytic of $\bm \alpha$ in a finite strip around $\mathbb T^d$ hence, by taking $|\bm \alpha|$ sufficiently close to the unit circle, $\|\rho_{\bm \alpha}(\bm r)\|<1$ and $\sum_{n=0}^\infty \big(\rho_{\bm \alpha}(\bm r)\big)^n$ is converges in norm together with its $\partial_\alpha$-derivatives. We conclude that $\rho_{\bm \alpha}\big ((1-r)^{-1}\big)$ is analytic in a strip around $\mathbb T^d$.\qed

\medskip $(\mathcal A \rtimes_\xi \mathbb Z^d)_{\mathrm{loc}}$ is a Fr\`{e}chet algebra and $(\mathcal A \rtimes_\xi \mathbb Z^d)_{\mathrm{loc}}$ belongs to the domain of any simple or higher derivation $\bm \partial^{\bm n}$. The stability under the holomorphic functional calculus ensures that the $\bm K_0$-groups of $(\mathcal A \rtimes_\xi \mathbb Z^d)_{\mathrm{loc}}$ and $\mathcal A \rtimes_\xi \mathbb Z^d$ coincide \cite{Connes:1994wk}.  Besides, as we shall see, the canonical Fredholm module associated to the crossed product $\mathcal A \rtimes_\xi \mathbb Z^d$ is automatically $(d+1)$-summable. These facts made us believe that $(\mathcal A \rtimes_\xi \mathbb Z^d)_{\mathrm{loc}}$ is the natural domain for the intrinsic Chern-Connes characters.

\section{THE REPRESENTATION}

In this Chapter we show that the $C^\ast$-algebra of adjointable operators over the standard Hilbert $\mathcal A$-module $\mathcal H_{\mathcal A}$ can serve as a natural imbedding algebra for $\mathcal A \rtimes_\xi \mathbb Z^d$. This will provide a natural connection between the $\bm K$-theories of $\mathcal A \rtimes_\xi \mathbb Z^d$ and of $\mathcal A$, as formulated within the framework of the standard Hilbert $\mathcal A$-module (see for example Chapter III in \cite{WeggeOlsenBook1993de}).

Viewing $\mathcal A$ as a right Hilbert $\mathcal A$-module, $\mathcal H_{\mathcal A}$ is defined as the tensor product $\mathcal A \otimes \mathcal H$ of Hilbert $C^*$-modules, with $\mathcal H$ being an ordinary separable Hilbert space \cite{BrucklerMC1999hf}.  The inner product on $\mathcal H_{\mathcal A}$ is $\langle a \otimes \phi| b \otimes \psi \rangle = \langle \phi | \psi \rangle a^\ast b$. The $C^\ast$-algebra of adjointable compact operators (\cite{KasparovJOT1980yt}, Def.~4) over $\mathcal H_{\mathcal A}$ is isomorphic to $\mathcal A \otimes \mathbb K$, where $\mathbb K$ is the algebra of compact operators over separable Hilbert spaces. The $C^\ast$-algebra of adjointable operators (\cite{KasparovJOT1980yt}, Def.~3) over $\mathcal H_{\mathcal A}$ is isomorphic to $\mathcal M(\mathcal A \otimes \mathbb K)$ (cf.~\cite{KasparovJOT1980yt}, Th.~1), the multiplier algebra or double centralizer \cite{BusbyTAMS1968jf} of $\mathcal A \otimes \mathbb K$. From a standard property of the tensor products and their multiplier algebras (\cite{WeggeOlsenBook1993de}, Corollary~T.6.3), we have a chain of algebra inclusions:
\begin{equation}\label{Inclusions}
\mathcal A \otimes \mathbb K \subset \mathcal A \otimes \mathbb B \subset \mathcal M(\mathcal A \otimes \mathbb K),
\end{equation}
where $\mathbb B$ is the algebra of bounded operators over separable Hilbert spaces and the tensor product  in $\mathcal A \otimes \mathbb B$ is with the spatial $C^\ast$-norm. This observation is important for us because the crossed-product algebra can be naturally imbedded in $\mathcal A \otimes \mathbb B$ and \ref{Inclusions} will provide an imbedding in $\mathcal M(\mathcal A \otimes \mathbb K)$, which is what we formally need.

In the present context, it is convenient to make the choice $\mathcal H = \ell^2(\mathbb Z^d)$, in which case the elements of $\mathcal H_{\mathcal A}$ can be uniquely expressed as:
\begin{equation}
(b_{\bm x})=\sum_{\bm x \in \mathbb Z^d} b_{\bm x} \otimes \delta_{\bm x} \in \mathcal A \otimes \ell^2(\mathbb Z^2),
\end{equation}
with $\{\delta_{\bm x}\}$ being the canonical orthonormal basis in $\ell^2(\mathbb Z^d)$. Similarly, the elements of $\mathcal A \otimes \mathbb K$ and $\mathcal A \otimes \mathbb B$ take the form:
\begin{equation}\label{Elements0}
\hat{\bm a} = \sum_{\bm x,\bm y} a_{\bm x \bm y} \otimes E_{\bm x, \bm y},
\end{equation}
with $E_{\bm x, \bm y}$ being the system of matrix-units over $\ell^2(\mathbb Z^d)$, $E_{\bm x, \bm y} \delta_{\bm z}= \delta_{\bm y \bm z} \delta_{\bm x}$. Let us also mention the ideal of finite-rank elements, algebraically generated by the rank-one operators:
\begin{equation}
\theta_{(a_{\bm x})}^{(b_{\bm x})} \big ( (c_{\bm x}) \big ) = (a_{\bm x}) \langle (b_{\bm x})|(c_{\bm x})\rangle, \ (a_{\bm x}), \ (b_{\bm x}), \ (c_{\bm x}) \in \mathcal H_{\mathcal A}.
\end{equation} 
In our settings:
\begin{equation}
\theta_{(a_{\bm x})}^{(b_{\bm x})} = \sum_{\bm x,\bm y \in \mathbb Z^d} a_{\bm x} b_{\bm y}^\ast \otimes E_{\bm x,\bm y}.
\end{equation} 
The ideal of finite-rank operators is dense in $\mathcal A \otimes \mathbb K$.

\begin{proposition} The following map
\begin{equation}\label{OurRep}
\mathcal A \rtimes_\xi \mathbb Z^d \ni \bm a \rightarrow \hat \pi (\bm a) = \sum_{\bm x,\bm q \in \mathbb Z^d} \xi_{\bm x} \Phi_{\bm q}(\bm a) \otimes E_{\bm x, \bm x - \bm q} \in \mathcal A \otimes \mathbb B,
\end{equation}
is well defined. It provides a faithful morphism of $C^\ast$-algebras and a faithful imbedding of $\mathcal A \rtimes_\xi \mathbb Z^d$ in $\mathcal M(\mathcal A \otimes \mathbb K)$.
\end{proposition}
\proof Our first task is to show that $\hat \pi$ is really into $\mathcal A \otimes \mathbb B$, {\emph i.e.} that $\hat \pi (\bm a)$ is a bounded operator when represented on some Hilbert space $\mathcal K \otimes \ell^2(\mathbb Z^d)$. For $\bm a \in \mathcal A \mathbb Z^d$, this can be accomplished by elementary means using a covariant representation of the dynamical system $(\mathcal A, \xi,\mathbb Z^d)$. Furthermore, for $\bm a$ and $\bm b$ from $\mathcal A \mathbb Z^d$, one can explicitly verify that (note the similarity between \ref{OurRep} and \ref{StandardRep}):
\begin{equation}
\hat \pi(\bm a) \hat \pi(\bm b) = \hat \pi (\bm a \bm b) \ \mathrm{and} \ \hat \pi(\bm a^\ast) = \hat \pi(\bm a)^\ast.
\end{equation}
The map is obviously faithful. At this point we established that $\hat \pi$ is a faithful $\ast$-representation of $\mathcal A \mathbb Z^d$ inside the $C^\ast$-algebra $\mathcal A \otimes \mathbb B$. By the very definition of $\mathcal A \rtimes_\xi \mathbb Z^d$, this representation must extend over the whole crossed-product.\qed

\section{THE CANONICAL TRACE AND ITS SPECIFIC CLASSES OF ELEMENTS}

The trace $\mathcal T_{\mathcal A}$ on $\mathcal A$ and the ordinary trace on $\mathbb B$ provide a trace $\widehat{\mathcal T}=\mathcal T_{\mathcal A} \otimes \mathrm{Tr}$ on $\mathcal A \otimes \mathbb B$. This trace can be characterized more directly as follows.

\begin{proposition}\label{Tr} Let $\mathcal A \otimes \mathbb B_+$ be the positive cone of $\mathcal A \otimes \mathbb B$. The functional:
\begin{equation}\label{Trace}
\mathcal A \otimes \mathbb B_+ \ni \hat{\bm a} \rightarrow \widehat{\mathcal T}\{\hat{\bm a}\} = \sum_{\bm x \in \mathbb Z^d} \mathcal T_{\mathcal A} \{a_{\bm x \bm x}\} \in [0,\infty],
\end{equation}
is a faithful, lower semicontinuous trace.
\end{proposition}

We recall the following standard consequences of $\widehat{\mathcal T}$ being a trace:\medskip 
\begin{itemize}
\item The set of trace-class elements:
\begin{equation}
\mathcal S_1 = \mathrm{span}\{\hat{\bm a} \in \mathcal A \otimes \mathbb B_+ \ | \ \widehat{\mathcal T}\{\hat{\bm a}\}<\infty \}
\end{equation}
is a two-sided ideal in $\mathcal A \otimes \mathbb B$.\medskip
\item The set of Hilbert-Schmidt elements:
\begin{equation}
\mathcal S_2 = \{\hat{\bm a} \in \mathcal A \otimes \mathbb B\ | \ \widehat{\mathcal T}\{\hat{\bm a}^\ast \hat{\bm a}\}<\infty \}
\end{equation}
is a two-sided ideal in $\mathcal A \otimes \mathbb B$ and $\mathcal S_1=\mathcal S_2 \cdot \mathcal S_2$.\medskip
\item The trace is cyclic:
\begin{equation}
\widehat{\mathcal T}\{\hat{\bm a} \hat{\bm b}\} = \widehat{\mathcal T}\{\hat{\bm b} \hat{\bm a}\}
\end{equation}
 for all $\hat{\bm a} \in \mathcal S_{\mathcal A}$ and $\hat{\bm b} \in \mathcal A \otimes \mathbb B$.\medskip
\item The trace is invariant to conjugation by unitaries:
\begin{equation}
\widehat{\mathcal T}\{\hat{\bm u}^\ast \hat{\bm a} \hat{\bm u} \} = \widehat{\mathcal T}\{\hat{\bm a}\}
\end{equation}
for all $\hat{\bm a} \in \mathcal S_{\mathcal A}$ and unitary $\hat{\bm u} \in \mathcal A \otimes \mathbb B$.
\end{itemize} 

The following characterization of the trace-class elements is essential for the generalized Fedosov principle elaborated in the next Chapter.

\begin{proposition}\label{TraceClass} $\mathcal S_1 = \mathcal A \otimes \mathbb S_1$, where $\mathbb S_1$ is the ideal of trace-class operators in $\mathbb B$. 
\end{proposition}
\proof Given that $\widehat{\mathcal T}$ is the extension of $\mathcal T_{\mathcal A} \otimes \mathrm{Tr}$, $\mathcal A \otimes \mathbb S_1$ automatically belongs to the domain of $\widehat{\mathcal T}$. We need to show the reverse inclusion, which follows if we can show that $\mathcal S_2 \subset \mathcal A \otimes \mathbb S_2$ where $\mathbb S_2$ is the ideal of Hilbert-Schmidt operators in $\mathbb B$. Since $\mathcal A$ is separable and $\mathcal T_{\mathcal A}$ is faithful and continuous, $\mathcal A$ can be completed to a separable Hilbert space when equipped with the inner product $\langle a | b \rangle = \mathcal T_{\mathcal A}(a^\ast b)$. Let $\{e_k\} \in \mathcal A$ be an orthonormal basis of this Hilbert space. Then any element in $\mathcal A \otimes \mathbb B$ can be uniquely represented as $\sum_{k} e_k \otimes T_k$. If such element belongs to $\mathcal S_2$, then 
\begin{equation}
\widehat{\mathcal T}\Big\{\Big(\sum_{k} e_k \otimes T_k  \Big)^\ast \sum_{k} e_k \otimes T_k \Big \} = \sum_k \mathrm{Tr}\{T_k^\ast T_k\} <\infty,
\end{equation}
hence  $\sum_{k} e_k \otimes T_k \in \mathcal A \otimes \mathbb S_2$.\qed

\medskip It will also be important to realize that the ideal of finite-rank operators belongs to the domain of the trace. Indeed:
\begin{equation}
\widehat{\mathcal T}\Big \{\theta_{(a_{\bm x})}^{(b_{\bm x})} \Big\} = \mathcal T_{\mathcal A} \big\{\langle (a_{\bm x}) | (b_{\bm x}) \rangle^\ast\big \}.
\end{equation}
Then it is obvious that any finite combination of $\theta$'s is in $\mathcal S_1$. 

\section{GENERALIZED FREDHOLM INDEX AND THE FEDOSOV FORMULA}

Within the framework of standard modules, $\bm K_0(\mathcal A)$ is defined as (\cite{WeggeOlsenBook1993de}, pg.~264):
\begin{equation}
\bm K_0(\mathcal A) =\{[\hat{\bm p}]-[\hat{\bm q}] \ | \ \hat{\bm p}, \ \hat{\bm q} \ \mbox{projectors in} \ \mathcal A \otimes \mathbb K \},
\end{equation}
where $[ \ ]$ indicates the classes under the Murray-von-Neumann equivalence relation in $\mathcal M(\mathcal A \otimes \mathbb K)$. We recall the index map for Hilbert $C^\ast$-modules, developed over a stretch of works by Kasparov \cite{KasparovJOT1980yt,KasparovMUI1981ui}, Mi\v{s}\v{c}enco and Fomenko \cite{FomenkoMUI1980js}, Pimsner, Popa, and Voiculescu \cite{PimsnerPopaVoiculescuJOT1980fj} and assembled in the final form by Mingo \cite{MingoTAMS1987fg}. Se also the pedagogical exposition in \cite{WeggeOlsenBook1993de} which calls this index the Mindex. We will do the same here. In short, one defines the generalized Calkin algebra as the corona $\mathcal M(\mathcal A \otimes \mathbb K)/\mathcal A \otimes \mathbb K$, and the class of generalized Fredholm elements as:
\begin{equation}
\mathcal F_{\mathcal A} = \{\hat{\bm f} \in \mathcal M(\mathcal A \otimes \mathbb K) \ | \ p(\hat{\bm f}) \ \mbox{invertible in} \ \mathcal M(\mathcal A \otimes \mathbb K)/\mathcal A \otimes \mathbb K\},
\end{equation}
where $p$ is the quotient map $\mathcal M(\mathcal A \otimes \mathbb K) \rightarrow \mathcal M(\mathcal A \otimes \mathbb K)/\mathcal A \otimes \mathbb K$. 

$\mathcal M(\mathcal A \otimes \mathbb K)$ is not a von-Neumann algebra, hence the polar decomposition cannot be assumed automatically.\footnote{According to an argument due to W.~J.~Phillips, the polar decomposition is equivalent to asking that $\hat{\bm f}$ has closed range (see~\cite{WeggeOlsenBook1993de}, Th.~15.3.8).} Nevertheless, any $\hat{\bm f} \in \mathcal F_{\mathcal A}$ admits a compact perturbation $\hat{\bm g}$ ({\it i.e.} $\hat{\bm f}-\hat{\bm g}\in \mathcal A \otimes \mathbb K$) which does accept a polar decomposition $\hat{\bm g}=\hat{\bm w} |\hat{\bm g}|$, with $\hat{\bm w}$ a partial isometry (see~\cite{MingoTAMS1987fg}, Proposition 1.7). This partial isometry defines two compact projectors: 
\begin{equation}
\mathrm{ker} \ \hat{\bm g} = 1-\hat{\bm w} \hat{\bm w}^\ast \in \mathcal A \otimes \mathbb K
\end{equation}
and
\begin{equation}
\mathrm{ker} \ \hat{\bm g}^\ast = 1 - \hat{\bm w}^\ast \hat{\bm w} \in \mathcal A \otimes \mathbb K,
\end{equation} 
which at their turn define an element of the $\bm K_0(\mathcal A)$-group:
\begin{equation}
[\mathrm{ker} \ \hat{\bm g}] - [\mathrm{ker} \ \hat{\bm g}^\ast] \in \bm K_0(\mathcal A).
\end{equation}
This element of $\bm K_0(\mathcal A)$ is insensitive of the compact perturbation $\hat{\bm g}$ used in the construction, hence it is truly determined by the Fredholm element $\hat{\bm f}$. Mingo's index map over $\mathcal H_{\mathcal A}$ is then defined as \cite{MingoTAMS1987fg}:
\begin{equation}
\mathrm{Mindex}\{\hat{\bm f}\} = [\mathrm{ker} \hat{\bm f}] - [\mathrm{ker} \hat{\bm f}^\ast] \in \bm K_0(\mathcal A),
\end{equation}
where it is understood that possible (irrelevant) compact perturbations are used to define the kernels. The Mindex is invariant to norm-continuous deformations of $\hat{\bm f}$ and two Fredholm operators have the same Mindex precisely when they are in the same path-component of $\mathcal F_{\mathcal A}$. The group of the homotopy classes of the generalized Fredholm elements characterizes completely the $\bm K_0$-group:

\begin{theorem}[\cite{MingoTAMS1987fg}]\label{MingoTh} $\bm K_0(\mathcal A) \simeq [\mathcal F_{\mathcal A}]$.
\end{theorem}

Our goal for this Chapter is to define a numerical index which is computationally more advantageous. We start with a statement which is standard for operators over ordinary Hilbert spaces: 

\begin{proposition} Any projector in $\mathcal A \otimes K$ is finite-rank.
\end{proposition}
\proof We will take advantage of the fact that the set of finite-rank elements is an ideal in $\mathcal M(\mathcal A \otimes \mathbb K)$ and that this ideal is dense in $\mathcal A \otimes \mathbb K$. Let
\begin{equation}\label{ApproxUnit}
\hat{\bm p}_N = \sum_{|\bm x| < N} 1 \otimes E_{\bm x, \bm x},
\end{equation}
be the standard approximation of identity for $\mathcal A \otimes \mathbb K$. Then, for any projector $\bm e \in \mathcal A \otimes \mathbb K$, $\hat{\bm e} \hat{\bm p}_N$ converges to $\hat{\bm e}$ in $\mathcal M(\mathcal A \otimes \mathbb K)$ as $N \rightarrow \infty$. Consequently, for $N$ large enough:
\begin{equation}
\hat{\bm s}=1-\hat{\bm e} (1-\hat{\bm p}_N)
\end{equation}
is invertible in $\mathcal M(\mathcal A \otimes \mathbb K)$ and $\hat{\bm e} \hat{\bm s} = \hat{\bm e}\hat{\bm p}_N$, hence a finite rank operator. Therefore $\hat{\bm e}=(\hat{\bm e}\hat{\bm p}_N) \hat{\bm s}^{-1}$ is finite-rank.\qed

\medskip This detail is important for our development because it shows that the compact projectors belong to the domain of the trace $\widehat{\mathcal T}$. As such, the pairing of $\widehat{\mathcal T}$ with $\bm K_0(\mathcal A)$ is straightforward:

\begin{proposition} The map:
\begin{equation}
\bm K_0(\mathcal A) \ni [\hat{\bm p}] -[\hat{\bm q}] \rightarrow \widehat{\mathcal T}\{\hat{\bm p}\} - \widehat{\mathcal T}\{\hat{\bm q}\} \in \mathbb R.
\end{equation}
is a group morphism.
\end{proposition}

\begin{definition} The generalized Fredholm index, which we call the $\widehat{\mathcal T}$-index, is defined as:
\begin{equation}
\mathcal F_{\mathcal A} \ni \hat{\bm f} \rightarrow \mathrm{Index}\{\hat{\bm f}\} = \widehat{\mathcal T}\{\mathrm{ker} \hat{\bm f}\} - \widehat{\mathcal T}\{\mathrm{ker} \hat{\bm f}^\ast\}.
\end{equation}
\end{definition}

The following characterization of the $\widehat{\mathcal T}$-index is a direct consequence of the properties of the Mindex mentioned earlier.

\begin{proposition} The $\widehat{\mathcal T}$-index defines a map $\widehat{\mathcal T}: \mathcal F_{\mathcal A} \rightarrow \mathbb R$ which is a locally constant homomorphism of semigroups (multiplicative for $\mathcal F_{\mathcal A}$ and additive for $\mathbb R$). 
\end{proposition}

As for the Breuer-Fredholm index \cite{CareyNCGPhysics}, the injectivity modulo path components present in \ref{MingoTh} may be lost and the information about $\bm K_0(\mathcal A)$ group, that can be extracted using the $\widehat{\mathcal T}$-index, can be limited. One recalls that this is not the case for the classic Fredholm index. However, there is an advantage for using both the $\widehat{\mathcal T}$-index and the Mindex, the former being algorithmically computable as the following generalization of the Fedosov's work \cite{FedosovTMMS1974fd} shows.

\begin{theorem}[Fedosov principle and formula]\label{FedosovTh} Let $\hat{\bm f} \in \mathcal M({\mathcal A \otimes \mathbb K})$ with $\|\hat{\bm f}\| \leq 1$. Then:
\begin{enumerate} 
\item If there exists a natural number $n$ such that:
\begin{equation}\label{Fedosov1}
(1 - \hat{\bm f}^\ast \hat{\bm f})^n \ \mathrm{and} \ (1 - \hat{\bm f} \hat{\bm f}^\ast)^n \in \mathcal S_1,
\end{equation}
then $\hat{\bm f} \in \mathcal F_{\mathcal A}$.
\item If (1) holds, then the  $\widehat{\mathcal T}$-index can be computed via:
\begin{equation}\label{Fedosov2}
\mathrm{Index}\{\hat{\bm f}\} = \widehat{\mathcal T}\{(1 - \hat{\bm f}^\ast \hat{\bm f})^n\} - \widehat{\mathcal T}\{(1 - \hat{\bm f} \hat{\bm f}^\ast)^n\}.
\end{equation}
\end{enumerate}
\end{theorem}

\proof ({\it i}) Assume \ref{Fedosov1} true. Then \ref{TraceClass} assures that $(1 - \hat{\bm f}^\ast \hat{\bm f})^n$ belongs to $\mathcal A \otimes \mathbb S_1$, hence compact. Taking $n'>n$ of the form $n'=2^k$, we can apply successively the square root on $(1 - \hat{\bm f}^\ast \hat{\bm f})^{n'}$ to access $1 - \hat{\bm f}^\ast \hat{\bm f}$. The latter is necessarily compact due to the following generic argument. Suppose $\hat{\bm a}\in \mathcal A \otimes \mathbb K$ is positive. Then its square root is defined as a positive element by the continuous functional calculus and:
\begin{equation}
\sqrt{\hat{\bm a}}= \lim_{\epsilon \searrow 0} \hat{\bm a}\big (\epsilon + \sqrt{\hat{\bm a}}\big )^{-1},
\end{equation}
where the limit is in $\mathcal M(\mathcal A \otimes \mathbb K)$. Since $\mathcal A \otimes \mathbb K$ is an ideal, the operators inside the limit belong to $\mathcal A \otimes \mathbb K$ and, since $\mathcal A \otimes \mathbb K$ is closed, the limit is in $\mathcal A \otimes \mathbb K$. Similarly, $1 - \hat{\bm f} \hat{\bm f}^\ast \in \mathcal A \otimes \mathbb K$. Therefore, $p(\hat{\bm f})$ is invertible in $\mathcal M(\mathcal A \otimes \mathbb K)/\mathcal A \otimes \mathbb K$.

({\it ii}) Assume for the beginning that $\hat{\bm f}$ accepts a polar decomposition $\hat{\bm f}=\hat{\bm w}|\hat{\bm f}|$, with $\hat{\bm w}= \hat{\bm w}\hat{\bm w}^\ast \hat{\bm w}$ and $\hat{\bm w}^\ast=\hat{\bm w}^\ast \hat{\bm w} \hat{\bm w}^\ast$. We write:
\begin{equation}
1 - \hat{\bm f} \hat{\bm f} ^\ast= 1 - \hat{\bm w}\hat{\bm w}^\ast +\hat{\bm w}(1 - |\hat{\bm f}|^2 ) \hat{\bm w}^\ast,
\end{equation}
and note that:
\begin{equation}
(1 - \hat{\bm w}\hat{\bm w}^\ast)\big  (\hat{\bm w}(1 - |\hat{\bm f}|^2) \hat{\bm w}^\ast \big )=0
\end{equation}
and
\begin{equation}
 \big (\hat{\bm w}(1 - |\hat{\bm f}|^2)\hat{\bm w}^\ast \big )  (1 - \hat{\bm w}\hat{\bm w}^\ast)=0.
\end{equation}
When combined with the elementary fact $\hat{\bm w}^\ast \hat{\bm w} |\hat{\bm f}| = |\hat{\bm f}|$, these lead to:
\begin{equation}
(1 - \hat{\bm f} \hat{\bm f}^\ast)^n = 1 - \hat{\bm w}\hat{\bm w}^\ast+\hat{\bm w}(1 - |\hat{\bm f}|^2 )^n \hat{\bm w}^\ast.
\end{equation}
Using the cyclic property of the trace:
\begin{equation}
\widehat{\mathcal T}\{\hat{\bm w}(1 - |\hat{\bm f}|^2 )^n \hat{\bm w}^\ast\}=\widehat{\mathcal T}\{\hat{\bm w}^\ast\hat{\bm w}(1 - |\hat{\bm f}|^2 )^n \}
\end{equation}
and, since $\hat{\bm w}^\ast \hat{\bm w} |\hat{\bm f}| = |\hat{\bm f}|$, it follows that:
\begin{equation}
\hat{\bm w}^\ast\hat{\bm w}(1 - |\hat{\bm f}|^2 )^n=\hat{\bm w}^\ast\hat{\bm w} - 1 + (1 - |\hat{\bm f}|^2 )^n.
\end{equation}
This together with \ref{Fedosov2} gives:
\begin{equation}
\widehat{\mathcal T}\{(1 - \hat{\bm f}^\ast \hat{\bm f})^n\} - \widehat{\mathcal T}\{(1 - \hat{\bm f} \hat{\bm f}^\ast)^n\}= \widehat{\mathcal T}\{1 - \hat{\bm w}^\ast \hat{\bm w}\} - \widehat{\mathcal T}\{1 - \hat{\bm w} \hat{\bm w}^\ast\}.
\end{equation}
The affirmation is then proven for the particular case when $\hat{\bm f}$ accepts a polar decomposition. 

For the generic case, note that we already established that $\hat{\bm f}$ is unitary modulo $\mathcal A \otimes \mathbb K$. Then we can apply Lemma~7.4 of \cite{PimsnerPopaVoiculescuJOT1980fj} which, as noted in \cite{MingoTAMS1987fg}, extends to the present context. This Lemma assures us that $\hat{\bm f}(1-\hat{\bm p}_N)$ accepts the polar decomposition, with $\hat{\bm p}_N$ defined in \ref{ApproxUnit} and $N$ large enough. Then $\hat{\bm f}(1-\hat{\bm p}_N)$ is a compact perturbation of $\hat{\bm f}$ with a polar decomposition, and since \ref{Fedosov1} still applies for $\hat{\bm f}(1-\hat{\bm p}_N)$:
\begin{equation}
\mathrm{Index}\{\hat{\bm f}\} = \widehat{\mathcal T}\{(1 - (1-\hat{\bm p}_N)\hat{\bm f}^\ast \hat{\bm f}(1-\hat{\bm p}_N))^n\} - \widehat{\mathcal T}\{(1 - \hat{\bm f} (1-\hat{\bm p}_N) \hat{\bm f}^\ast)^n\}.
\end{equation}
Using the cyclic property of the trace, one can check directly that all the terms containing $\hat{\bm p}_N$ cancel identically, hence this last equation can be reduced to \ref{Fedosov2}.\qed

\medskip The above statements can be generalized in the following way. Suppose there are projections $\hat{\bm p}$ and $\hat{\bm p}'$ such that $\hat{\bm p}' \hat{\bm f} = \hat{\bm f} \hat{\bm p}= \hat{\bm f}$ and $\hat{\bm p}'=\hat{\bm s} \hat{\bm p} \hat{\bm s}$ with $\hat{\bm s}$ a symmetry (\emph{i.e.} $\hat{\bm s}$ self-adjoint and $\hat{\bm s}^2=1$). Let:
\begin{equation}
\tilde{\bm f} = (1-\hat{\bm p})\hat{\bm s}+\hat{\bm f}.
\end{equation}
By applying the generalized Fedosov principle on $\tilde{\bm f}$, we learn that $\tilde{\bm f}\in \mathcal F_{\mathcal A}$ provided:
\begin{equation}\label{Fedosov3}
( \hat{\bm p}'- \hat{\bm f}^\ast \hat{\bm f})^n \ \mathrm{and} \ (\hat{\bm p} - \hat{\bm f} \hat{\bm f}^\ast)^n \in \mathcal S_1.
\end{equation}
Furthermore, if \ref{Fedosov3} holds, then we can apply the generalized Fedosov formula on $\tilde{\bm f}$, which gives:
\begin{equation}\label{Fedosov4}
\mathrm{Index}\{\tilde{\bm f}\} = \widehat{\mathcal T}\{(\hat{\bm p}' - \hat{\bm f}^\ast \hat{\bm f})^n\} - \widehat{\mathcal T}\{(\hat{\bm p} - \hat{\bm f} \hat{\bm f}^\ast)^n\}.
\end{equation}
The above index will be understood as the index of $\hat{\bm f}$ as an operator between the right $\mathcal A$-modules $\hat{\bm p}\mathcal H_{\mathcal A}$ and $\hat{\bm p}'\mathcal H_{\mathcal A}$.

\section{THE INTRINSIC FREDHOLM MODULE}

We extend $\hat \pi$ to a representation $\hat \pi_\gamma =\mathrm{id} \otimes \hat \pi$ on $Cl_d \otimes \mathcal M(\mathcal A \otimes \mathbb K)$, where $Cl_{d}$ is the even Clifford algebra with generators: 
\begin{equation}
\gamma_i \gamma_j + \gamma_j \gamma_i = 2\delta_{ij} \  i, j=1,\ldots,d. 
\end{equation}
We denote by $\mathrm{Tr}_\gamma$ the standard normalized trace on $Cl_d$. Kasparov's stabilization theorem (\cite{KasparovJOT1980yt}, Th.~2) assures us that we can replace $\mathcal A$ by $Cl_d \otimes \mathcal A$ without changing the $\bm K$-theory and, by examining the arguments in the previous Chapter, one can convince himself that replacing $\mathcal T_{\mathcal A}$ by $\mathrm{Tr}_\gamma \otimes \mathcal T_{\mathcal A}$ will fully accommodate the new setup. 

This extension enables us to define the Dirac element:
\begin{equation}
\widehat{\bm d}_{\bm x_0} = \sum_{\bm x \in \mathbb Z^d} \bm \gamma \cdot(\bm x +\bm x_0) \otimes 1 \otimes E_{\bm x,\bm x},
\end{equation}
where $\bm x \cdot \bm y$ denotes the Euclidean scalar product $\bm x \cdot \bm y = x_1 y_1 + \ldots + x_d y_d$. Note that $\widehat{\bm d}_{\bm x_0} $ is an unbounded operator over $\mathcal H_{\mathcal A}$ but its phase $\widetriangle{\bm d}_{\bm x_0}=\mathrm{sign}(\widehat{\bm d}_{\bm x_0})$ is part of the imbedding algebra:
\begin{equation}
\widetriangle{\bm d}_{\bm x_0} = \sum_{\bm x \in \mathbb Z^d} \bm \gamma \cdot \widetriangle{\bm x +\bm x_0} \otimes 1 \otimes E_{\bm x,\bm x} \in Cl_d \otimes \mathcal A \otimes \mathbb B.
\end{equation}
Throughout, we will use the notation $\widetriangle{\bm x}=\bm x/|\bm x|$.

\begin{remark} The shift $\bm x_0$ is allowed to take values in the cube $[0,1]^d$. When $\bm x_0 \notin \mathbb Z^d$, the shift is useful because it provides an un-ambiguous phase. However, the considerations leading to the inclusion of the shift $\bm x_0$ go beyond that. In particular, the average over $\bm x_0 \in [0,1]^d$ is absolutely essential for the key identity \ref{CentralId}. When $\bm x_0 \in \mathbb Z^d$ then $ \bm \gamma \cdot \widetriangle{\bm x + \bm x_0} $ needs to be modified at one point. This can be done in various ways and is not an important detail. 
\end{remark}

One can verify directly that $(\widetriangle{\bm d}_{\bm x_0})^\ast = \widetriangle{\bm d}_{\bm x_0}$ and $(\widetriangle{\bm d}_{\bm x_0})^2=1$. Also:
\begin{equation}
\gamma = \gamma_0 \otimes 1 \otimes 1, \ \gamma_0 = -\imath^{\frac{d}{2}} \gamma_1 \gamma_2 \ldots \gamma_d,
\end{equation}
provides a grading that has the right properties: 
\begin{equation}\label{P87}
\widetriangle{\bm d}_{\bm x_0} \gamma = - \gamma \widetriangle{\bm d}_{\bm x_0}, \ \mathrm{and} \ [\hat \pi_\gamma(\bm a),\gamma]=0 \ \forall \ \bm a\in \mathcal A \rtimes_\xi \mathbb Z^d.
\end{equation} 
All these lead to the following definition.

\begin{definition}\label{GFredholmModule} The family of intrinsic Fredholm modules for $\mathcal A \rtimes_\xi \mathbb Z^d$ is defined as:
\begin{equation}
(Cl_d \otimes \mathcal M(\mathcal A \otimes \mathbb K), \hat \pi_\gamma, \widetriangle{\bm d}_{\bm x_0},\gamma)_{\bm x_0 \in [0,1]^d}.
\end{equation}
\end{definition}

\begin{definition}\label{SummaDef} The family of Fredholm modules is $n$-summable over a sub-algebra of $\mathcal A \rtimes_\xi \mathbb Z^d$ if:
\begin{equation}
\prod_{i=1}^n \big [\widetriangle{\bm d}_{\bm x_0},\hat \pi_\gamma(\bm a_i) \big ] \in Cl_d \otimes \mathcal S_1,
\end{equation}
for any $\bm a_i$ in that sub-algebra. 
\end{definition}

\begin{theorem}\label{Summability} The intrinsic family of Fredholm modules is $n$-summable over $(\mathcal A \rtimes_\xi \mathbb Z^d)_{\mathrm{loc}}$, for any $n \geq d+1$.
\end{theorem}

\proof Let $\bm a_i \in (\mathcal A \rtimes_\xi \mathbb Z^d)_{\mathrm{loc}}$. We will first show that:
\begin{equation}\label{S1}
\Big\| \Big ( \prod_{i=1}^k \big [\widetriangle{\bm d}_{\bm x_0},\hat \pi_\gamma(\bm a_i) \big ] \Big )_{\bm x \bm x}\Big \| \leq A  |\bm x + \bm x_0|^{-k}.
\end{equation}
Consequently, $\prod_{i=1}^k \big [\widetriangle{\bm d}_{\bm x_0},\hat \pi_\gamma(\bm a_i) \big ]$ belongs to $\mathcal S_2$ for all $k >d/2$ and the main affirmation follows for $n \geq d+2$. Additional work will be needed to cover the case $n=d+1$. Returning to \ref{S1}, it is useful to write the commutators explicitly:
\begin{equation}
\big [\widetriangle{\bm d}_{\bm x_0},\hat \pi_\gamma(\bm a_i)\big ] = \sum_{\bm x,\bm y \in \mathbb Z^d} \bm \gamma \cdot (\widetriangle{ \bm x + \bm x_0} - \widetriangle{\bm y + \bm x_0}) \otimes \xi_x \Phi_{\bm x - \bm y}(\bm a_i) \otimes E_{\bm x, \bm y}.
\end{equation}
Then
\begin{align}\label{PP32}
& \Big (\prod_{i=1}^k \big [\widetriangle{\bm d}_{\bm x_0},\hat \pi_\gamma(\bm a_i) \big ]\Big )_{\bm x \bm x}=(\mathrm{id} \otimes \xi_{\bm x})\sum_{\bm x_i \in \mathbb Z^d} \delta_{\bm x_1,\bm 0} \delta_{\bm x_{k+1},0}   \\
 & \nonumber \ \ \ \prod_{i=1}^k \bm \gamma \cdot(\widetriangle{ \bm x_i + \bm x + \bm x_0} - \widetriangle{\bm x_{i+1} + \bm x + \bm x_0}) \otimes \xi_{\bm x_i} \Phi_{\bm x_i - \bm x_{i+1}}(\bm a_i),
\end{align}
and
\begin{align}
& \Big \| \Big (\prod_{i=1}^k \big [\widetriangle{\bm d}_{\bm x_0},\hat \pi_\gamma(\bm a_i) \big ]\Big )_{\bm x \bm x} \Big \| \leq A \sum_{\bm x_i \in \mathbb Z^d} \delta_{\bm x_1,\bm 0} \delta_{\bm x_{k+1},0}   \\
& \nonumber \ \ \ \prod_{i=1}^k \big |\widetriangle{ \bm x_i + \bm x + \bm x_0} - \widetriangle{\bm x_{i+1} + \bm x + \bm x_0} \big | \| \Phi_{\bm x_i - \bm x_{i+1}}(\bm a_i)\|.
\end{align} 
Throughout, all uninteresting constants will be denoted by $A$. Due to the asymptotic behavior:
\begin{equation}\label{Asy}
\widetriangle{ \bm y + \bm x } - \widetriangle{\bm y' + \bm x } \sim |\bm x|^{-1} \Big ( \bm y - \bm y' +\big (\widetriangle{\bm x}\cdot (\bm y - \bm y')\big ) \widetriangle{\bm x} \Big ),
\end{equation}
as $|\bm x|\rightarrow \infty$, the supremum 
\begin{equation}
S(\bm y,\bm y')=\sup \big \{|\bm x| \big |\widetriangle{ \bm y + \bm x } - \widetriangle{\bm y' + \bm x} \big |, \bm x \in \mathbb R^d  \big \}
\end{equation}
is finite. $S(\bm y,\bm y')$ has the scaling property:
\begin{equation}
S(s \bm y, s \bm y') = s S(\bm y,\bm y'),
\end{equation}
hence, by taking $s = (|\bm y|+|\bm y'|)^{-1}$, we obtain the upper bound:
\begin{equation}\label{S5}
S(\bm y,\bm y') \leq (|\bm y|+|\bm y'|) \sup\{S(\bm x,\bm x'), \ |\bm x|+|\bm x'| = 1\}.
\end{equation}
The conclusion is:
\begin{align}
 & \Big \| \Big (\prod_{i=1}^k \big [\widetriangle{\bm d}_{\bm x_0},\hat \pi_\gamma(\bm a_i) \big ]\Big )_{\bm x \bm x} \Big \| \leq | \bm x + \bm x_0|^{-k} \\
& \nonumber \ \ \ A\sum_{\bm x_i \in \mathbb Z^d} \delta_{\bm x_1,\bm 0} \delta_{\bm x_{k+1},0} 
\prod_{i=1}^k (|\bm x_i|+|\bm x_{i+1}|) \alpha^{|\bm x_i - \bm x_{i+1}|},
\end{align} 
with $\alpha <1$. The remaining sum is finite.

Let us now consider the element:
\begin{equation}
\hat{\bm e} = \sum_{\bm x \in \mathbb Z^d} |\bm x|^{-\frac{1}{4}} \ 1 \otimes 1 \otimes E_{\bm x, \bm x} \in Cl_d \otimes \mathcal A \otimes \mathbb B.
\end{equation}
Given \ref{S1}, we have:
\begin{equation}
\Big ( \prod_{i=1}^{\nicefrac{d}{2}} \big [\widetriangle{\bm d}_{\bm x_0},\hat \pi_\gamma(\bm a_i) \big ] \Big )\hat{\bm e} \in \mathcal S_2, \ \forall \ \bm a_i \in (\mathcal A \rtimes_\xi \mathbb Z^d)_{\mathrm{loc}}
\end{equation}
and, by conjugation, also:
\begin{equation}\label{S57}
\hat{\bm e} \Big ( \prod_{i=1}^{\nicefrac{d}{2}} \big [\widetriangle{\bm d}_{\bm x_0},\hat \pi_\gamma(\bm a_i) \big ] \Big )\in \mathcal S_2,  \ \forall \ \bm a_i \in (\mathcal A \rtimes_\xi \mathbb Z^d)_{\mathrm{loc}}.
\end{equation}
Consider the sum:
\begin{equation}
\hat{\bm b} = \sum_{\bm x,\bm q \in \mathbb Z^d} |\bm x|^{\frac{1}{4}} \langle \bm \gamma, \widetriangle{ \bm x + \bm x_0} - \widetriangle{\bm x -\bm q + \bm x_0}\rangle |\bm x -\bm q|^{\frac{1}{4}} \otimes \xi_x \Phi_{\bm q}(\bm a_{\frac{d}{2}+1}) \otimes E_{\bm x, \bm x - \bm q},
\end{equation}
with $\bm a_{\frac{d}{2}+1} \in (\mathcal A \rtimes_\xi \mathbb Z^d)_{\mathrm{loc}}$. Using \ref{S5}, we have:
\begin{equation}
\|b_{\bm x,\bm x - \bm q} \| \leq A |\bm x|^{-\frac{1}{2}} \alpha^{|\bm q|}, \ \alpha < 1,
\end{equation}
hence $\hat{\bm b} \in Cl_d \otimes \mathcal A \otimes \mathbb B$. Since $\mathcal S_2$ is an ideal:
\begin{equation}
\Big ( \prod_{i=1}^{\nicefrac{d}{2}} \big [\widetriangle{\bm d}_{\bm x_0},\hat \pi_\gamma(\bm a_i) \big ] \Big )\hat{\bm e} \hat{\bm b} \in \mathcal S_2, \ \forall \ \bm a_i \in (\mathcal A \rtimes_\xi \mathbb Z^d)_{\mathrm{loc}}
\end{equation}
and combining with \ref{S57}:
\begin{equation}
\Big ( \prod_{i=1}^{\nicefrac{d}{2}} \big [\widetriangle{\bm d}_{\bm x_0},\hat \pi_\gamma(\bm a_i) \big ] \Big )\hat{\bm e} \hat{\bm b}\hat{\bm e} \Big ( \prod_{i=\nicefrac{d}{2}+1}^{d+1} \big [\widetriangle{\bm d}_{\bm x_0},\hat \pi_\gamma(\bm a_i) \big ] \Big ) \in \mathcal S_1.
\end{equation}
A direct computation will show that:
\begin{equation}
\hat{\bm e} \hat{\bm b}\hat{\bm e} = \big [\widetriangle{\bm d}_{\bm x_0},\hat \pi_\gamma(\bm a_{\frac{d}{2}+1}) \big ],
\end{equation}
and the main affirmation follows, with $n=d+1$ this time.\qed

\section{THE INTRINSIC CHERN-CONNES CHARACTER AND ITS LOCAL FORMULA}

\begin{definition}\label{ChernConnesCh} The intrinsic Chern-Connes character is defined by the cohomology class of following $d+1$-cyclic cocycle:
\begin{equation}
\tau_d(\bm a_0,\ldots,\bm a_d) = \nicefrac{1}{2}\int\limits_{[0,1]^d}d\bm x_0 \ \widehat{\mathcal T}\Big \{\gamma \widetriangle{\bm d}_{\bm x_0} \prod_{i=0}^d \big [\widetriangle{\bm d}_{\bm x_0}, \hat \pi_\gamma(\bm a_i)\big ] \Big \},
\end{equation}
defined over $(\mathcal A \rtimes_\xi \mathbb Z^d)_{\mathrm{loc}}$.
\end{definition}

\begin{remark} We recall that the classic Chern-Connes characters are defined using a representation and the trace on an ordinary Hilbert space \cite{Connes:1994wk}. Same algebra is required to show that $\tau_d$ above is a cocycle.
\end{remark}

\begin{theorem}[The $KK$-map]\label{KKMap} Let $\bm p \in (\mathcal A \rtimes_\xi \mathbb Z^d)_{\mathrm{loc}}$ be a projector and let 
$$\hat{\pi}_\gamma^\pm = \nicefrac{1}{2}(1 \pm \gamma)\hat{\pi}_\gamma$$ 
be the decomposition of the representation $\hat{\pi}$ with respect to the grading $\gamma$. Then the element $\hat \pi^-_\gamma(\bm p) \hat{\bm d}_{\bm x_0}\hat \pi_\gamma^+ (\bm p) \in \mathcal M(\mathcal A \otimes \mathbb K)$ is Fredholm, hence it provides a map:
\begin{equation}
\bm K_0(\mathcal A \rtimes_\xi \mathbb Z^d) \ni [\bm p]\rightarrow \mathrm{Mindex}\{\hat \pi^-_\gamma(\bm p) \widetriangle{\bm d}_{\bm x_0}\hat \pi_\gamma^+ (\bm p)\} \in \bm K(\mathcal A).
\end{equation}
\end{theorem}
\proof Let us list a few useful identities:
\begin{equation}\label{Id07}
\widetriangle{\bm d}_{\bm x_0} \big [ \widetriangle{\bm d}_{\bm x_0},\hat \pi_\gamma(\bm p) \big ]=- \big [ \widetriangle{\bm d}_{\bm x_0},\hat \pi_\gamma(\bm p) \big ]\widetriangle{\bm d}_{\bm x_0},
\end{equation}
\begin{equation}\label{Id08}
\hat \pi_\gamma(\bm p) \big [ \widetriangle{\bm d}_{\bm x_0},\hat \pi_\gamma(\bm p) \big ]=- \big [ \widetriangle{\bm d}_{\bm x_0},\hat \pi_\gamma(\bm p) \big ]\hat \pi_\gamma(\bm p)
\end{equation}
and
\begin{equation}\label{Id11}
\hat \pi_\gamma(\bm p) \big [ \widetriangle{\bm d}_{\bm x_0},\hat \pi_\gamma(\bm p) \big ] ^2 = \hat \pi_\gamma(\bm p) \widetriangle{\bm d}_{\bm x_0} \hat \pi_\gamma(\bm p) \widetriangle{\bm d}_{\bm x_0} \hat \pi_\gamma(\bm p) -\hat \pi_\gamma(\bm p).    
\end{equation}
Now, if we take $\hat{\bm f}=\hat \pi^-_\gamma(\bm p) \hat{\bm d}_{\bm x_0}\hat \pi_\gamma^+ (\bm p)$, then obviously we have the projectors $\hat \pi^\pm_\gamma(\bm p)$ such that $\hat \pi^+_\gamma(\bm p) \hat{\bm f} = \hat{\bm f} \hat \pi^-_\gamma(\bm p)=\hat{\bm f}$ and $(\gamma_1\otimes 1\otimes 1) \hat \pi^\pm_\gamma(\bm p) (\gamma_1\otimes 1\otimes 1) = \hat \pi^\mp_\gamma(\bm p)$. As such, we are allowed to use the Fedosov principle in the form \ref{Fedosov3}. From identity \ref{Id11}:
\begin{align}\label{KK32}
&(\hat \pi^-_\gamma(\bm p) - \hat{\bm f} \hat{\bm f}^\ast)^n = - \hat \pi_\gamma^-(\bm p) \big [\widetriangle{\bm d}_{\bm x_0},\hat \pi_\gamma(\bm p) \big ]^{2n} \\
\label{KK33} &(\hat \pi^+_\gamma(\bm p) - \hat{\bm f}^\ast \hat{\bm f})^n = - \hat \pi_\gamma^+(\bm p) \big [ \widetriangle{\bm d}_{\bm x_0},\hat \pi_\gamma(\bm p) \big ]^{2n},
\end{align}
and, given the summability result \ref{Summability}, the Fedosov principle applies if we take $2n>d+1$.\qed

\begin{theorem}[The index theorem] \label{IndexTh}
\begin{equation}\label{Index0}
\tau_d(\bm p,\ldots,\bm p)= \mathrm{Index}\big \{ \hat \pi^-_\gamma(\bm p) \hat{\bm d}_{\bm x_0}\hat \pi_\gamma^+ (\bm p) \big \}.
\end{equation}
The $\widehat{\mathcal T}$-index on the right is independent of $\bm x_0$.
\end{theorem}
\proof Given \ref{Asy}, $\hat{\bm d}_{\bm x_0} - \hat{\bm d}_{\bm x_0'}$ is compact, hence the index is independent of $\bm x_0$. We take $n=\frac{d}{2}+1$ and, aided by \ref{KK32} and \ref{KK33}, we apply the Fedosov-formula \ref{Fedosov4}:
\begin{equation}
\mathrm{Index}\big \{ \hat \pi^+_\gamma(\bm p) \hat{\bm d}_{\bm x_0}\hat \pi_\gamma^- (\bm p) \big \} =-\widehat{\mathcal T} \Big \{ \gamma \hat \pi_\gamma(\bm p) \big [ \widetriangle{\bm d}_{\bm x_0},\hat \pi_\gamma(\bm p) \big ]^{d+2} \Big \}.
\end{equation}
Using \ref{Id11}, the righthand side becomes:
\begin{equation}
-\widehat{\mathcal T} \Big \{ \gamma  \Big ( \hat \pi_\gamma(\bm p) \widetriangle{\bm d}_{\bm x_0} \hat \pi_\gamma(\bm p) \widetriangle{\bm d}_{\bm x_0} \hat \pi_\gamma(\bm p) -\hat \pi_\gamma(\bm p) \Big ) \big [ \widetriangle{\bm d}_{\bm x_0},\hat \pi_\gamma(\bm p) \big ]^{d} \Big \},
\end{equation}
which can be rewritten in two different ways:
\begin{equation}
\widehat{\mathcal T} \Big \{ \gamma \widetriangle{\bm d}_{\bm x_0} \hat \pi_\gamma(\bm p) \big [ \widetriangle{\bm d}_{\bm x_0},\hat \pi_\gamma(\bm p) \big ]^{d+1} \Big \},
\end{equation}
or
\begin{equation}
\widehat{\mathcal T} \Big \{ \gamma \widetriangle{\bm d}_{\bm x_0} \big [ \widetriangle{\bm d}_{\bm x_0},\hat \pi_\gamma(\bm p) \big ]^{d+1} \hat \pi_\gamma(\bm p) \Big \}.
\end{equation}
Eq.~\ref{Index0} then follows from the following elementary identity:
\begin{equation}
\pi_\gamma(\bm p) \big [  \widetriangle{\bm d}_{\bm x_0},\hat \pi_\gamma(\bm p) \big ]^{d+1}  +\big [  \widetriangle{\bm d}_{\bm x_0},\hat \pi_\gamma(\bm p) \big ]^{d+1}\hat \pi_\gamma(\bm p) =\big [  \widetriangle{\bm d}_{\bm x_0},\hat \pi_\gamma(\bm p) \big ]^{d+1},
\end{equation}
property \ref{P87} and an average over $\bm x_0$. \qed

\begin{theorem}[The local formula] \label{LocalFormulaTh}
The Chern-Connes cocycle accepts the local formula:
\begin{equation}
\tau_d(\bm a_0,\ldots,\bm a_d)= \Lambda_d \sum_{\rho \in S_d} (-1)^\rho \mathcal T\Big \{\bm a_0  \prod_{i=1}^d \partial_{\rho_i} \bm a_i \Big \},
\end{equation}
where $(\bm \partial,\mathcal T)$ is the non-commutative calculus over $\mathcal A \rtimes_\xi \mathbb Z^d$ and $\Lambda_d$ was defined in \ref{L}.
\end{theorem}

\proof Opening the first commutator in the product of definition \ref{ChernConnesCh} transforms its righthand side into:
\begin{align}
\nicefrac{1}{2} \ \sum_{\bm x \in \mathbb Z^d} \mathrm{Tr}_\gamma \otimes \mathcal T_{\mathcal A} \Big \{ \Big (\gamma \big (\hat{\pi}_\gamma(\bm a_0)  - \widetriangle{\bm d}_{\bm x_0} \hat{\pi}_\gamma(\bm a_0) \widetriangle{\bm d}_{\bm x_0} \big ) \prod_{i=1}^{d} [\widetriangle{\bm d}_{\bm x_0},\hat{\pi}_\gamma(\bm a_i)]\Big )_{\bm x \bm x}\Big \}.
\end{align}
Given \ref{Id07} and since $d$ is even, we can swap $\widetriangle{\bm d}_{\bm x_0}$ and the product, and using the cyclic property of $\mathrm{Tr}_\gamma$ we arrive at:
\begin{align}
\tau_d(\bm a_0,\ldots,\bm a_d) =\int\limits_{[0,1]^d} d\bm x_0\sum_{\bm x \in \mathbb Z^d} \mathrm{Tr}_\gamma \otimes \mathcal T_{\mathcal A} \Big \{ \Big (\gamma \hat{\pi}_\gamma(\bm a_0) \prod_{i=1}^{d} [\widetriangle{\bm d}_{\bm x_0},\hat{\pi}_\gamma(\bm a_i)]\Big )_{\bm x \bm x}\Big \}.
\end{align}
Now, using \ref{PP32}:
\begin{align}
& \Big (\gamma \hat{\pi}_\gamma(\bm a_0)  \prod_{i=1}^d \big [\widetriangle{\bm d}_{\bm x_0},\hat \pi_\gamma(\bm a_i) \big ]\Big )_{\bm x \bm x}  \\
 & \nonumber \ \ \ = (\mathrm{id} \otimes \xi_{\bm x})\sum_{\bm x_i \in \mathbb Z^d} \delta_{\bm x_{d+1},0}  \  \gamma_0 \otimes \Phi_{-\bm x_1}(\bm a_0)\\
& \nonumber \ \ \ \ \ \ \prod_{i=1}^{d} \bm \gamma \cdot (\widetriangle{ \bm x_i + \bm x + \bm x_0} - \widetriangle{\bm x_{i+1} + \bm x + \bm x_0}) \otimes \xi_{\bm x_i} \Phi_{\bm x_i - \bm x_{i+1}}(\bm a_i).
\end{align}
Using the invariance of $\mathcal T_{\mathcal A}$ with respect to $\xi$-automorphisms and by combining the integration over $\bm x_0$ with the summation over $\bm x$, we conclude:
\begin{align}
& \tau_d(\bm a_0,\ldots,\bm a_d)  = \sum_{\bm x_i \in \mathbb Z^d} \delta_{\bm x_{d+1},0} \\
& \nonumber \ \ \ \int_{\mathbb R^d} d\bm x    \ \mathrm{Tr}_\gamma \Big \{ \gamma_0 \prod_{i=1}^d \bm \gamma \cdot (\widetriangle{\bm x +\bm x_i} -\widetriangle{\bm x +\bm x_{i+1}}) \Big \}  \\
& \nonumber \ \ \ \ \ \ \mathcal T_{\mathcal A} \Big \{\Phi_{- \bm x_1}(\bm a_0) \prod_{i=1}^d \xi_{\bm x_i} \Phi_{\bm x_i - \bm x_{i+1}}(\bm a_i) \Big \}.
\end{align}
At this point we use the identity discovered in \cite{ProdanJPA2013hg}:
\begin{align}\label{CentralId}
&\delta_{\bm x_{d+1},0}\int\limits_{\mathbb R^d} d{\bm x} \ \mathrm{Tr}_\gamma \Big \{ \gamma_0 \prod_{i=1}^d \bm \gamma \cdot (\widetriangle{\bm x+ \bm x_i} -\widetriangle{\bm x + \bm x_{i+1}}) \Big \} = \tilde \Lambda_d \sum_{\rho \in S_d} (-1)^\rho \prod_{i=1}^d (\bm x_i)_{\rho_i},
\end{align}
with $\tilde \Lambda_d = -\frac{(2\pi)^\frac{d}{2}}{\imath^\frac{d}{2} (d/2)!}$, to continue:
\begin{align}
\tau_d(\bm a_0,\ldots,\bm a_d) & = \tilde \Lambda_d \sum_{\rho \in S_d} (-1)^\rho \sum_{\bm x_i \in \mathbb Z^d}  \mathcal T_{\mathcal A} \Big \{\Phi_{- \bm x_1}(\bm a_0) \prod_{i=1}^d (\bm x_i)_{\rho_i} \xi_{\bm x_i} \Phi_{\bm x_i - \bm x_{i+1}}(\bm a_i) \Big \}.
\end{align}
Due to the anti-symmetrizer factor, we can replace  $(\bm x_i)_{\rho_i}$ by $ (\bm x_i)_{\rho_i}- (\bm x_{i+1})_{\rho_i}$ and then $\big ( (\bm x_i)_{\rho_i}- (\bm x_{i+1})_{\rho_i} \big ) \Phi_{\bm x_i - \bm x_{i+1}}(\bm a_i)$ by $  -\imath \Phi_{\bm x_i - \bm x_{i+1}}(\partial_{\rho_i}\bm a_i)$. After these substitutions we can recognize that:
\begin{align}
\tau_d(\bm a_0,\ldots,\bm a_d) & =  \Lambda_d \sum_{\rho \in S_d} (-1)^\rho   \mathcal T_{\mathcal A} \Big \{\Phi_{0}\Big(\bm a_0 \prod_{i=1}^d  \partial_{\rho_i}\bm a_i \Big) \Big \}
\end{align}
and the affirmation follows.\qed

\section{AN APPLICATION}\label{Applications}

One of the most successful applications of the noncommutative geometry in condensed matter is the solution to the Integer Quantum Hall Effect (IQHE) \cite{BELLISSARD:1994xj}, which explained the quantization and homotopy invariance of the Hall conductance of a 2-dimensional electron gas subjected to a perpendicular magnetic field and strong disorder.  It was Haldane \cite{HaldanePRL1988rh} who first realized that certain materials can exhibit all the characteristics of the standard IQHE even in the absence of an external magnetic field. This type of materials are now called Chern insulators. Presently, there exists an entire classificaltion table of presumably \emph{all} topological insulating phases of matter \cite{SchnyderPRB2008qy,kitaev:22,RyuNJP2010tq}. When examining this table (see for example Table III in \cite{RyuNJP2010tq}), one notices that there are no Chern insulators in 3 space-dimensions. The physical argument is that, although the transport coefficients do show some topological characteristics in the absence or at weak disorder, those features will presumably disappear in the regime of strong disorder. Translating, this means that, although $\bm K_0(C(\Omega) \rtimes_\xi \mathbb Z^3) = \mathbb Z^4$ [$C(\Omega) \rtimes_\xi \mathbb Z^3$ is the algebra of physical observables, see below], the cocycles that can be paired with $K_0$ are of degree lower than the space dimension, hence they are not stable once the spectral gap of the Hamiltonian closes due to strong disorder. From a mathematical point of view, as far as we know, this problem is completely open. The issue is definitely interesting, being also tied to the quantized Hall-Effect in 3-dimensions, theoretically proposed quite some time ago \cite{KohmotoPB1993hf} but whose decisive experimental confirmation is still to come.

The quantum dynamics of electrons in homogeneous 3-dimensional materials is generated by self-adjoint operators on $\ell^2(\mathbb Z^3,\mathbb C^Q)$ of the form \cite{BellissardLN2003bv}:
\begin{equation}\label{DisorderH}
(H\psi)_{\bm x}=\sum_{\bm q \in \mathbb Z^3} \big (1+\lambda_{\bm q}(\xi_{\bm x}\omega) \big ) A_{\bm q}  \psi_{\bm x-\bm q},
\end{equation}
where $\lambda_{\bm q}: \Omega \rightarrow \mathbb C$ are functions over a classical dynamical system $(\Omega,\xi)$ and $A_{\bm q}$'s are $Q\times Q$ ordinary matrices called hopping matrices. Typically, $\lambda_{\bm q}$'s are much smaller than 1, hence the particular writing in \ref{DisorderH}. One will recognize in \ref{DisorderH} the standard representation $\pi_\omega h$ (see \ref{StandardRep}) of:
\begin{equation}
h=\sum_{\bm q} \phi_{\bm q} \cdot \bm q \in C(\Omega) \otimes M_Q(\mathbb C) \rtimes_\xi \mathbb Z^3,
\end{equation}
with:
\begin{equation}
\phi_{\bm q}(\omega) = \big ( 1+\lambda_{\bm q}(\omega)\big ) \otimes A_{\bm q}.
\end{equation}

For disordered crystals, the prototypical $\Omega$ is a Tychonoff space:
\begin{equation}
\Omega = {\mathcal I}^{\mathbb Z^3}, \ \Omega \ni \omega= \{\omega_{\bm x}\}_{\bm x \in \mathbb Z},
\end{equation}
and the action of $\mathbb Z^3$ is provided by the shift: $\xi_{\bm y}(\omega)=\{\omega^{\bm q}_{\bm x+\bm y}\}$. $\mathcal I$ is endowed with a probability measure $dP_{\mathcal I}$ and the product measure $d\mathbb P(\omega) = \prod_{\bm x \in \mathbb Z^3} dP_{\mathcal I}(\omega_{\bm x})$ provides an ergodic and $\xi$-invariant probability measure over $\Omega$. At its turn, $d \mathbb P(\omega)$ provides a natural trace over $C(\Omega)$, hence over $C(\Omega)\otimes M_Q(\mathbb C)$. For topological insulators, one can distinguish two different disorder regimes \cite{ProdanJPhysA2011xk}: 1) The weak disorder regime, where $\lambda_{\bm q}$'s are small and some spectral gaps of $H$ remain open. The Fermi level $\epsilon_F$ is fixed in the middle of such a spectral gap. 2) The strong disorder regime, where $\lambda_{\bm q}$'s are large and all the spectral gaps of $H$ are closed but there are still regions of pure-point spectrum. $\epsilon_F$ is fixed in middle of such a region.

Although we are not yet in the position to resolve the strong disorder regime, the machinery developed in this work is quite relevant here, as it generates an index formula for the transport coefficients. The linear conductivity tensor $\sigma$ is defined by the relation $J_i =\sigma_{ij}E_j$, where $\bm J$ is the electron current set in motion by a weak electric field $\bm E$. The off-diagonal components of $\sigma$ were computed in \cite{KohmotoPB1993hf}. In the language of crossed product algebras, it takes a form similar to that of the 1-st noncommutative Chern number:
\begin{equation}\label{Sigma}
\sigma_{12} =  2\pi \imath \sum_{\rho \in S_2} (-1)^\rho \mathcal T_3\Big \{\bm p  \prod_{i=1}^2 \partial_{\rho_i} \bm p \Big \},
\end{equation}
in some adjusted physical units, but note the mismatch between the degree of the cycle and the dimension of the crossed product. Above, $\mathcal T_3$ is the trace \ref{Trace0} on $C(\Omega)\otimes M_Q(\mathbb C) \rtimes_\xi \mathbb Z^3$ and $\bm p=\chi_{(-\infty,\epsilon_F]}(h)$ is the projector onto the spectrum of $h$ up to $\epsilon_F$. The existence of a spectral gap at $\epsilon_F$, \emph{i.e} the weak disorder regime, ensures that $\bm p$ belongs to the $C^\ast$-algebra $C(\Omega)\otimes M_Q(\mathbb C) \rtimes_\xi \mathbb Z^3$. 

One can already see that the righthand side of \ref{Sigma} is a pairing between a cyclic-cocycle and a projector, hence a homotopy invariant \cite{Connes:1994wk}. We can now complete with an index formula which will tell us where this pairing takes place. We want to point out that the methods of \cite{ProdanJPA2013hg} cannot be used here because the cocycle is not trace-class when represented on a Hilbert space! Within the new framework, we identify $\mathcal A = C(\Omega)\otimes M_Q(\mathbb C) \rtimes_{\xi_3} \mathbb Z$ and define $\mathcal T_{\mathcal A}$ via the Fourier calculus on $C(\Omega)\otimes M_Q(\mathbb C) \rtimes_{\xi_3} \mathbb Z$ (cf.~\ref{Trace0}), and similarly for the trace $\mathcal T$ over $\mathcal A \rtimes_{\xi_{12}} \mathbb Z^2$. Then $C(\Omega)\otimes M_Q(\mathbb C)\rtimes_\xi \mathbb Z^3 = \mathcal A \rtimes_{\xi_{12}} \mathbb Z^2$ and $\mathcal T_3 = \mathcal T$. Hence:
\begin{equation}
\sigma_{12} =  2 \pi \imath \sum_{\rho \in S_2} (-1)^\rho \mathcal T\Big \{\bm a_0  \prod_{i=1}^2 \partial_{\rho_i} \bm p \Big \}=\mathrm{Index}\{\hat{\pi}_\gamma^-(\bm p) \widetriangle{d}_{\bm x_0} \hat{\pi}_\gamma^+(\bm p)\}.
\end{equation}
Consequently, $\sigma_{12}$ takes values in $\widehat{\mathcal T}\big (K_0(C(\Omega) \rtimes_{\xi_3} \mathbb Z)\big )$. This provides a prediction about the pattern of possible experimental values for $\sigma_{12}$.

\bibliographystyle{plain}
\bibliography{../../../TopologicalInsulators}

\end{document}